\documentclass[11pt]{amsart}
\usepackage{amssymb,amsmath,epsfig,graphics}
\theoremstyle{plain}
\newtheorem{thrm}{Theorem}[section]
\newtheorem{lemma}[thrm]{Lemma}
\newtheorem{prop}[thrm]{Proposition}
\newtheorem{cor}[thrm]{Corollary}
\newtheorem{rmrk}[thrm]{Remark}
\newtheorem{dfn}[thrm]{Definition}

\numberwithin{equation}{section} \numberwithin{figure}{section}

\begin{document}

\title[Geometric second derivative estimates in Carnot groups and convexity]
{Geometric second derivative estimates in Carnot groups and
convexity}

\author{Nicola Garofalo}
\address{Department of Mathematics\\Purdue University \\
West Lafayette, IN 47907} \email[Nicola
Garofalo]{garofalo@math.purdue.edu}

%\address{Dipartimento di Metodi e Modelli Matematici per le Scienze Applicate\\ Universit\`a di Padova\\ 35131 Padova, Italy}
%\email[Nicola Garofalo]{garofalo@dmsa.unipd.it}

\thanks{Supported in part by NSF Grant DMS-07010001}

\maketitle

% begin top matter
% ********************* macroes needed for this paper ************************
\newcommand{\SL}{\mathcal L^{1,p}(\Om)}
\newcommand{\Lp}{L^p(\Omega)}
\newcommand{\CO}{C^\infty_0(\Omega)}
\newcommand{\Rn}{\mathbb R^n}
\newcommand{\Rm}{\mathbb R^m}
\newcommand{\R}{\mathbb R}
\newcommand{\Om}{\Omega}
\newcommand{\Hn}{\mathbb H^n}
\newcommand{\HH}{\mathbb H^1}
\newcommand{\eps}{\epsilon}
\newcommand{\BVX}{BV_H(\Omega)}
\newcommand{\IO}{\int_\Omega}
\newcommand{\bG}{\mathbb{G}}
\newcommand{\bg}{\mathfrak g}
\newcommand{\p}{\partial}
\newcommand{\Xnu}{\overset{\rightarrow}{ H_\nu}}
\newcommand{\nuX}{\boldsymbol{\nu}^H}
\newcommand{\Up}{\bN^H}
\newcommand{\n}{\boldsymbol \nu}
\newcommand{\sigmau}{\boldsymbol{\sigma}^u_H}
\newcommand{\di}{\nabla_{i}^{H,\p \Om}}
\newcommand{\del}{\delta_H}
\newcommand{\nui}{\nu_{H,i}}
\newcommand{\nuj}{\nu_{H,j}}
\newcommand{\dej}{\delta_{H,j}}
\newcommand{\cx}{\boldsymbol{c}_S}
\newcommand{\sx}{\sigma_H}
\newcommand{\lx}{\mathcal L_H}
\newcommand{\rad}{\overline u}
\newcommand{\nH}{\nabla_H}
\newcommand{\uk}{u(\frac{|x|^2}{4},|y'|)}
\newcommand{\nh}{\nabla^H}
\newcommand{\uc}{u_{,ij}}
\newcommand{\sij}{\sum_{i,j=1}^m}
\newcommand{\sul}{\Delta_H}
\newcommand{\nor}{\boldsymbol \nu}
\newcommand{\bN}{\boldsymbol N}
\newcommand{\pb}{\overline p}
\newcommand{\qb}{\overline q}
\newcommand{\ob}{\overline \omega}
\newcommand{\mS}{\mathcal S}
\newcommand{\om}{\omega}
\newcommand{\delXY}{\nabla^{H,\mS}_X Y}
\newcommand{\oX}{\overline X}
\newcommand{\oY}{\overline Y}

% ****************************************************************************

\begin{abstract}
We prove some new a priori estimates for $H_2$-convex functions
which are zero on the boundary of a bounded smooth domain $\Om$ in a
Carnot group $\bG$. Such estimates are global and are geometric in
nature as they involve the horizontal mean curvature $\mathcal H$ of
$\p \Om$. As a consequence of our bounds we show that if $\bG$ has
step two, then for any smooth $H_2$-convex function in $\Om \subset
\bG$ vanishing on $\p \Om$ one has
\begin{equation}\label{commHab}
\sij \int_\Om ([X_i,X_j]u)^2\ dg \ \leq\  \frac{4}{3} \int_{\p \Om}
\mathcal H\ |\nabla_H u|^2\ d\sigma_H\ . \notag
\end{equation}
\end{abstract}

\section{\textbf{Introduction and statement of the results}}\label{S:I}

\vskip 0.2in

In this paper we study some a priori estimates of geometric type for
functions vanishing on the boundary of a smooth bounded open set in
a Carnot group. The estimates that we obtain are of interest in
connection with the study of the geometric notion of convexity in
Carnot groups recently introduced in \cite{DGN}, see also
\cite{LMS}, \cite{GM}, \cite{GT}, \cite{DGNT}, \cite{Wa1},
\cite{Wa2}, \cite{Ma}. They also play a central role in establishing
a priori bounds in $L^2$ for the (horizontal) second derivatives of
solutions of non-variational operators with rough coefficients. In
this respect, a global version of our results (with sharp constant)
for compactly supported functions in the Heisenberg group $\Hn$ has
been recently established in the very interesting recent work of
Domokos and Manfredi, see Lemma 1.1 in \cite{DM}, using the deeper
spectral decomposition of Strichartz \cite{Str}. An alternative
proof of the estimates in \cite{DM}, which however does  not produce
the best constant, would be to use the subelliptic estimates, see
\cite{Ko} and \cite{Ho}.

Our results are intimately connected to those in \cite{DGNT}, except
that the approach in that paper was based on monotonicity formulas
for a certain fully nonlinear subelliptic operator, rather than
geometric inequalities as those in the present paper. One
interesting and novel aspect here is that the relevant estimates
depend in an explicit way on a new geometric object, the so called
\emph{horizontal mean curvature} of the boundary of the ground
domain. Whether our method is capable of producing estimates with
sharp constants or not presently remains an open question. Our
approach is based on some delicate (but otherwise fairly elementary)
integration by parts formulas which are combined with a
sub-Riemannian Bochner type identity. The latter is inspired by the
classical one from Riemannian geometry which states that for a
Riemannian manifold $M$ with Levi-Civita connection $\nabla $, one
has for $u\in C^3(M)$
\begin{equation}\label{bochner}
\Delta (|\nabla u|^2)\ =\ 2\ ||\nabla^2 u||^2\ +\ 2\ <\nabla
u,\nabla (\Delta u)>\ + \ 2\ Ric(\nabla u,\nabla u)\ ,
\end{equation}
where $Ric(\cdot,\cdot)$ represents the Ricci tensor. A beautiful
generalization to CR manifolds of \eqref{bochner} was found by
Greenleaf in \cite{Gre}. Another tool that we use in the proof of
Theorem \ref{T:est} below is a sub-Riemannian Rellich identity
discovered in \cite{GV}.

In the classical case a basic a priori estimate of the elliptic
theory reads
\begin{equation}\label{be} ||u||_{W^{2,2}(\Om)}\ \leq\ C\
\left(||Lu||_{L^2(\Om)}\ +\ ||u||_{L^2(\Om)}\right)\ ,
\end{equation}
where $\Om\subset \Rn$, and $L$ is a second order uniformly elliptic
operator. One of the main tools in the obtainment of \eqref{be} is
the following geometric a priori inequality
\begin{equation}\label{et} \int_\Om ||\nabla^2 u||^2 \ dx\ \leq\ C\
\bigg(\int_\Om (\Delta u)^2\ dx\ +\ \int_\Om u^2\ dx\bigg)\ ,
\end{equation} valid for any $u\in C^2(\Om)$, with $u = 0$ on $\p \Om$,
where $\Om\subset \Rn$ is a bounded open set which is piecewise
$C^1$ and whose principal curvatures are bounded, see \cite{LU},
Lemma 8.1 on p.175. Here $\nabla^2 u$ denotes the Hessian matrix of
$u$, and $||\nabla^2 u||$ indicates its Hilbert-Schmidt norm. The
constant $C>0$ depends on various parameters, among which
appropriate bounds on the principal curvatures of $\p \Om$. The
prototype of estimates such as \eqref{be} and \eqref{et} first
appeared in two dimension in the pioneering works of S. N. Bernstein
\cite{Be1}, \cite{Be2}. Several years later Kadlec \cite{Ka} first
obtained a higher dimensional version of \eqref{et} for convex
domains. One should also see the works of Ladyzenskaya and Uraltseva
\cite{LU}, Talenti \cite{Ta}, Grisvard \cite{Gr}, Lewis \cite{Le}.

We consider a Carnot group $\bG$ with Lie algebra $\bg$ (for the
relevant definitions and properties see Section \ref{S:prelim}). We
assume throughout that $\bG$ is endowed with a left-invariant
Riemannian metric with respect to which the Lie algebra generating
left-invariant vector fields $X_1,...,X_m$ defined in \eqref{vf}
below constitute an orthonormal basis of the horizontal subbundle
$H\bG$. If $\Om\subset \bG$ is an open set and $k\in \mathbb N$, we
denote by $\Gamma^k(\Om)$ the Folland-Stein class of functions $u\in
C(\bG)$ such that for every $1\leq s\leq k$ one has
$X_{j_1}...X_{j_s}u \in C(\bG)$, where $j_{\ell}\in \{1,...,m\}$ for
$1\leq \ell \leq s$. Given a function $u\in \Gamma^2(\bG)$ we let
$\nabla^H u = \sum_{i=1}^m X_i u X_i$ denote the horizontal gradient
of $u$. If $\zeta = \sum_{i=1}^m \zeta_i X_i \in
\Gamma^1(\bG,H\bG)$, then we let $div_H \zeta = \sum_{i=1}^m X_i
\zeta_i$. The horizontal Laplacian of $u\in \Gamma^2(\bG)$ is then
given by
\[ \Delta_H u\ =\ div_H \ \nabla^H u\ =\ \sum_{i=1}^m X_i^2 u\ .
\]

The symmetrized Hessian of $u$ is the $m\times m$ matrix defined by
\[ \nabla^2_H u\ =\ [u_{,ij}]\ ,\quad\quad\quad\text{where}\quad
u_{,ij}\ =\ \frac{X_iX_j u + X_j X_i u}{2}\ , \] and clearly
$\Delta_H u = trace\ \nabla_H^2 u$. The matrix $\nabla_H^2 u$ plays
a central role in the study of convexity in Carnot groups. It was in
fact proved in \cite{DGN}, \cite{LMS} that a function $u\in
\Gamma^2$ is $H$-convex (see definition \eqref{convex0} below) if
and only if $\nabla_H^2u\geq 0$.

Our first result is an integral identity  which connects an
interesting fully nonlinear subelliptic operator to the geometry of
the ground domain through the $H$-mean curvature $\mathcal H$ of its
boundary (for the latter notion, see Definition \ref{D:HMC} below).
In what follows we denote by $d\sigma_H$ the $H$-perimeter measure
on $\p \Om$, see Definition \ref{D:permeas} and also \cite{DGN2}.

\begin{thrm}\label{T:grisvard}
Let $\bG$ be a Carnot group and consider a $C^2$ bounded open set
$\Om\subset \bG$. Let $u\in \Gamma^3(\overline \Om) \cap
C^2(\overline \Om)$  with $u\leq 0$ in $\Om$ and $u = 0$ on $\p
\Om$. One has
\begin{align}\label{1}
& \int_\Om \bigg\{(\Delta_H u)^2 -  ||\nabla^2_H u||^2\bigg\} dg\
 +\ \frac{3}{4}\ \sum_{i,j=1}^m \int_\Om ([X_i , X_j] u)^2\ dg\\
&  +\ \sij \int_\Om X_iu\ [[X_i,X_j],X_j]u\ dg\ =\ \int_{\p \Om}
\mathcal H\ |\nh u|^2 d\sigma_H\ . \notag
\end{align}
When the step of $\bG$ is two, then the third integral in the
left-hand side of \eqref{1} vanishes. If instead $u\in
C^\infty_0(\Om)$, then one obtains
\begin{align}\label{2}
\int_\Om ||\nabla^2_H u||^2\ dg\ & =\ \int_\Om (\Delta_H u)^2\ dg\
+\ \frac{3}{4} \sij \int_\Om ([X_i,X_j]u)^2\ dg
\\
&  +\ \sij \int_\Om X_iu \big[[X_i,X_j],X_j\big]u\ dg\ . \notag
\end{align}
\end{thrm}

It is interesting to state explicitly Theorem \ref{T:grisvard} in
the special, yet important situation, when $\bG = \Hn$, the
Heisenberg group, with group law
\[ g \circ g' = (x,y,t) \circ (x',y',t') = (x + x',y+y',t+t' +
\frac{1}{2}(x\cdot y' - x' \cdot y))\ ,
\]
and the left-invariant basis for the Lie algebra
\begin{equation}\label{vf}
X_j = \frac{\partial}{\partial x_j} - \frac{y_j}{2}\
\frac{\partial}{\partial t},\quad\quad  X_{n+j} =
\frac{\partial}{\partial y_{j}} + \frac{x_j}{2}\
\frac{\partial}{\partial t},\quad j=1,...,n,\ \ \ \ T\ =\ \frac{\p
}{\p t}\ .
\end{equation}

We note that the vector fields $X_1,...,X_{2n}$ satisfy the
commutation relations
\[ [X_j, X_{n+k}]\ =\ \delta_{jk}\ T\ ,\ \ \  j, k = 1,...,n\ ,
\]
and therefore they generate the Lie algebra of $\Hn$.

\begin{cor}\label{C:grisvardHn}
Consider a $C^2$ bounded open set $\Om\subset \Hn$. Let $u\in
\Gamma^3(\overline \Om) \cap C^2(\overline \Om)$ with $u\leq 0$ in
$\Om$ and $u = 0$ on $\p \Om$. One has
\begin{align}\label{1Hn}
& \int_\Om (\Delta_H u)^2\ dg\ -\ \int_\Om ||\nabla^2_H u||^2\ dg +\
\frac{3}{2}\ n \int_\Om (Tu)^2\ dg
\\
& =\ \int_{\p \Om} \mathcal H\ |\nabla_H u|^2\ d\sigma_H\ . \notag
\end{align}
As a consequence, if $\p \Om$ has $H$-mean curvature $\mathcal H
\geq 0$, then
\begin{equation}\label{2Hn}
 \int_\Om ||\nabla^2_H u||^2\ dg\ \leq \ \int_\Om (\Delta_H u)^2\
dg\ +\ \frac{3}{2}\ n  \int_\Om (Tu)^2\ dg \ .
\end{equation}
If instead $u\in C^\infty_0(\Om)$, then regardless of the sign of
$\mathcal H$ one obtains
\begin{equation}\label{3Hn}
\int_\Om ||\nabla^2_H u||^2\ dg\  =\ \int_\Om (\Delta_H u)^2\ dg\ +\
\frac{3}{2}\ n \int_\Om (Tu)^2\ dg \ .
\end{equation}
When $n=1$, then \eqref{1Hn} becomes
\begin{equation}\label{4}
 \int_\Om det(\nabla^2_H u)\ dg\ +\ \frac{3}{4} \int_\Om (Tu)^2\ =\ \frac{1}{2}\ \int_{\p \Om} \mathcal H\ |\nabla_H u|^2\ d\sigma_H\
.
\end{equation}
\end{cor}

With some appropriate modifications, Theorems \ref{T:grisvard} and
Corollary \ref{C:grisvardHn} are still valid if one removes the
assumption that $u\leq 0$ in $\Om$. The following result easily
follows by keeping track of the various terms appearing in the proof
of Theorem \ref{T:grisvard}.

\begin{thrm}\label{T:grisvardstep2bis}
Let $\bG$ be a Carnot group of step $r=2$, and consider a $C^2$
bounded open set $\Om\subset \bG$. Let $u\in \Gamma^3(\overline \Om)
\cap C^2(\overline \Om)$  with $u = 0$ on $\p \Om$. One has
\begin{align}\label{1bis}
&  \left|\int_\Om (\Delta_H u)^2\ dg\ -\ \int_\Om ||\nabla^2_H
u||^2\ dg\ +\ \frac{3}{4} \sij \int_\Om ([X_i,X_j]u)^2\ dg\right|
\\
&  \leq\  \int_{\p \Om} |\mathcal H|\ |\nabla_H u|^2\ d\sigma_H\ .
\notag
\end{align}
\end{thrm}

We next obtain some basic consequences of Theorem \ref{T:grisvard}
when the function $u$ is $H_2$-\emph{convex}. With this hypothesis
we are able to bound the $L^2$ norm of the commutators $[X_i,X_j]u$
in terms of a weighted integral of the $H$-mean curvature of $\p
\Om$.

To introduce the relevant notions we recall that for $r=1,...,m$,
the $r$-th elementary symmetric function is defined by
\begin{equation}\label{sf0}
S_r(x)\ =\ \underset{i_1 < ... < i_r}{\sum} x_{i_1}\ ...\ x_{i_r}\
,\quad\quad\quad\quad 1 \leq r \leq m\ .
\end{equation}

When $r>1$ we can use such functions to form the fully nonlinear
differential operators
\begin{equation}\label{Fr}
\mathcal F_r[u]\ =\ S_r(\lambda_1(u), ... , \lambda_m(u))\ ,
\end{equation}
where $\lambda_1(u),...,\lambda_m(u)$ denote the eigenvalues of the
symmetrized Hessian of $u$. One easily recognizes that
\begin{equation}\label{F1}
\mathcal F_1[u]\ =\ S_1(\lambda)\ =\ trace(\nabla^2_H u)\ =\ \sul u\
=\ \sum_{i=1}^m u_{,ii}\ \ \quad\quad \text{(horizontal Laplacian)}\
,
\end{equation}
\begin{equation}\label{F2}
\mathcal F_2[u]\ =\ S_2(\lambda)\ =\ \sum_{i<j} \left(u_{,ii}\
u_{,jj}\ -\ u_{,ij}^2\right)\ =\ \frac{1}{2}\ \left\{(\sul u)^2\ -
||\nabla^2_H u||^2\right\}\ ,
\end{equation}
\[
.\ \ \ \ \ .\ \ \ \ \ .
\]
\[
.\ \ \ \ \ .\ \ \ \ \ .
\]
\[
.\ \ \ \ \ .\ \ \ \ \ .
\]
\begin{equation}\label{Fm}
\mathcal F_m[u]\ =\ S_m(\lambda)\ =\ \det\ \nabla^2_H(u)\ \quad\quad
\text{(horizontal Monge-Amp\`ere)}\ .
\end{equation}

\begin{dfn}\label{D:sf}
For $r = 1,...,m$, a function $u\in \Gamma^2(\bG)$ is called
$H_r$-convex, if $\mathcal F_k(u) \geq 0$ for $k=1,...,r$.
\end{dfn}

For these notions and for related results we refer the reader to the
paper \cite{DGNT}.

\begin{rmrk}\label{R:convex}
We observe that $H_1$-convex functions correspond to subharmonic
functions, i.e., $\sul u \geq 0$, whereas a function $u$ is
$H_2$-convex if $\sul u\geq 0$ and $(\sul u)^2 - ||\nabla^2_H u||^2
\geq 0$. We recall the following geometric notion of convexity
introduced in \cite{DGN}. One should also see \cite{LMS} where a
notion of convexity in the viscosity sense was set forth. These two
notions have been recently shown to be equivalent, see \cite{Wa1},
\cite{Wa2}, \cite{Ri}. A function $u:\bG \to \R$ is called
$H$-\emph{convex} if given any point $g\in \bG$ and $0 \leq \lambda
\leq 1$, the following inequality holds
\begin{equation}\label{convex0}
u(g \delta_\lambda (g^{-1} g'))\ \leq\ (1 - \lambda) u(g)\ +\
\lambda u(g')\ ,\quad\quad\quad\text{for every}\quad g'\in H_g\ ,
\end{equation}
where $H_g$ indicates the horizontal plane through $g\in \bG$. In
\eqref{convex0} we have denoted by $\delta_\lambda : \bG \to \bG$
the anisotropic dilations on $\bG$. The point $g \delta_\lambda
(g^{-1} g')$ denotes the twisted convex combination of $g$ and $g'$
based at $g$. According to Theorem 5.12 in \cite{DGN}, a function
$u$ is $H_m$-convex according to Definition \ref{D:sf} if and only
if $u$ is $H$-convex.
\end{rmrk}

\begin{thrm}\label{T:commH}
Let $\bG$ be a Carnot group of step $r=2$, and consider a $C^2$
bounded open set $\Om\subset \bG$. Let $u\in \Gamma^3(\overline \Om)
\cap C^2(\overline \Om)$ be $H_2$-convex in $\Om$ with $u = 0$ on
$\p \Om$. One has
\begin{equation}\label{commH1}
\sij \int_\Om ([X_i,X_j]u)^2\ dg \ \leq\  \frac{4}{3} \int_{\p \Om}
\mathcal H\ |\nabla_H u|^2\ d\sigma_H\ . \notag
\end{equation}
In particular, if $\bG = \mathbb H^n$, one obtains
\begin{equation}\label{commH2}
\int_\Om (Tu)^2\ \leq\ \frac{2}{3n}\ \int_{\p \Om} \mathcal H\
|\nabla_H u|^2\ d\sigma_H\ .
\end{equation}
\end{thrm}

We note that, according to Lemma \ref{L:convlevelset}, under the
hypothesis of Theorem \ref{T:commH} we must have $\mathcal H \geq 0$
on $\p \Om \setminus \Sigma$, where $\Sigma$ denotes the
characteristic set of $\p \Om$, see definition \eqref{csaf} below.
Since thanks to results of Balogh \cite{Ba} and Magnani \cite{Ma} we
know that $\sigma_H(\Sigma) = 0$, we conclude that $\mathcal H \geq
0$ $\sigma_H$-a.e. on $\p \Om$.

One should compare the sharp geometric bounds in Theorem
\ref{T:commH} with the following non-geometric local a priori bound
established in \cite{DGNT}, see also \cite{GM}, \cite{GT} and
\cite{GM2}.

\begin{thrm}\label{T:oscillation}
Consider a bounded open set $\Om$ in a group of step two $\bG$. Let
$u\in \Gamma^3(\Om)$ be a $H_2$-convex function. For any $D\subset
\subset D'\subset \subset \Om$ we have for some constant $C>0$
depending on $\bG,\Om$, $D'$, and $D$
\[
\sum_{i,j=1}^m\ \int_D  \left([X_i,X_j]u\right)^2\ dg\ \leq\ C\
\left(\underset{D'}{osc}\ u\right)^{2}\ .
\]
\end{thrm}

The next theorem provides a basic global a priori bound for the
$L^2$ norms of the commutators of an $H_2$-convex function vanishing
on the boundary under the assumption that the ground domain be
starlike (in a weak sense) and that the horizontal mean curvature of
its boundary be bounded. In the starlikeness assumption in
\eqref{starlike} below, the vector field $\mathcal Z$ indicates the
infinitesimal generator of the non-isotropic group dilations. For
its definition see \eqref{Z} below.

\begin{thrm}\label{T:est}
Let $\bG$ be a Carnot group of step $r=2$, and let $\Om\subset \bG$
be a $C^2$ bounded open set such that for some $M , \alpha>0$,
\begin{equation}\label{boundedmc}
\underset{\p \Om}{\sup}\ |\mathcal H|\ \leq\ M\ ,
\end{equation}
and, with $W$ as in \eqref{ps} below, suppose that
\begin{equation}\label{starlike}
\underset{\p \Om}{\inf}\ <\mathcal Z,\n>\ \geq\ \alpha\ W\ .
\end{equation}
There exists a constant $C(\bG,\Om, M , \alpha)>0$ such that for
$u\in \Gamma^3(\overline \Om) \cap C^2(\overline \Om)$ which is
$H_2$-convex in $\Om$, and satisfies $u = 0$ on $\p \Om$, one has
\begin{equation}\label{kadlec}
 \sij  \int_{\Om} ([X_i,X_j] u)^2\ dg \
 \leq \ C\ \int_\Om\ (\Delta_H u)^2\ dg\  .
\end{equation}
\end{thrm}

The proof of Theorem \ref{T:est} is accomplished by combining
Theorem \ref{T:commH} with a sub-Riemannian Rellich identity
discovered in \cite{GV}, see Theorem \ref{T:Rellich} and Corollaries
\ref{C:Rellich} and \ref{C:Rellich2}. These results allow to
establish Lemma \ref{L:rellich3} which is instrumental in
controlling the commutator term in \eqref{kadlec} in Theorem
\ref{T:est} solely in terms of the $L^2$ norm of the horizontal
Laplacian. We notice explicitly that, since by \eqref{csaf} at the
characteristic points of $\p \Om$ the \emph{angle function} $W$
vanishes, condition \eqref{starlike} is a weak starlikeness
assumption with respect to the non-isotropic group dilations
\eqref{dil} below. We also emphasize that in a Carnot group of
arbitrary step a basic family of domains satisfying the hypothesis
\eqref{boundedmc} and \eqref{starlike} in Theorem \ref{T:est} is
represented by the gauge pseudo-balls centered at the group
identity, see Proposition \ref{P:gaugeball}.

We mention at this point that after this paper was submitted we
learnt from J. Manfredi of his interesting preprint \cite{CM} joint
with Chanillo in which, using the above mentioned Bochner identity
due to Greenleaf \cite{Gre}, the authors obtain a priori estimates
connected to those in Theorem \ref{T:grisvard}, but for strictly
pseudo-convex CR hypersurfaces $M^{2n+1}\subset \mathbb C^{n+1}$.
There are however two essential differences between our work and
\cite{CM}, and neither of these papers is contained in the other. On
the one hand there is the obvious fact that there exists in nature a
plentiful supply of Carnot groups which are not CR hypersurfaces.
The second distinction has to do with the different goals of the
papers. To explain this point we mention that, as it is well-known,
see e.g. \cite{S}, the Heisenberg group $\Hn$ with its flat CR
structure, via its identification with the boundary of the Siegel
upper half-space, is the basic prototype of a CR hypersurface.
However, even in this specialized context the results in \cite{CM}
do not contain ours since we are primarily concerned with global
geometric estimates connecting the fully nonlinear operators
introduced in \cite{DGNT} to the horizontal mean curvature of a
relatively compact sub-domain of the group itself. On the other
hand, in \cite{CM} in the non-compact case the authors work
exclusively with $C^\infty_0(M^{2n+1})$ functions and the geometry
of the boundary plays no role for them. It appears that the ultimate
goal in \cite{CM} is achieving the sharp constant in their
commutator estimates since this allows them to generalize to the CR
setting the above cited Cordes type results in \cite{DM}.

In closing we briefly describe the organization of the paper.
Section \ref{S:prelim} is devoted to recalling some basic facts
about Carnot groups and the notion of horizontal (or
sub-Riemannnian) mean curvature of a hypersurface in such a group.
In section \ref{S:bochner} we prove the Bochner identity in
Proposition \ref{P:B}. In section \ref{S:2nd} we prove our main
results: Theorems \ref{T:grisvard}, \ref{T:commH} and \ref{T:est}.
Finally, in Proposition \ref{P:gaugeball} we show that in any Carnot
group the gauge pseudo-balls satisfy the hypothesis of Theorem
\ref{T:est}.

\noindent \textbf{Acknowledgment:} The results in this paper were
presented in the special session \emph{Subelliptic PDE's and
Sub-Riemannian Spaces} at the AMS Fall Southeastern Section Meeting,
University of Arkansas, November 3-4, 2006. The author thanks the
organizers L. Capogna, S. Pauls and J. Tyson for their gracious
invitation.

\section{\textbf{Preliminaries}}\label{S:prelim}

 Consider a Carnot group $\bG$ of step $r$. This is a simply connected
  Lie group whose Lie algebra $\bg$ is graded and $r$-nilpotent. This means that there exists vector sub-spaces $V_1,...,V_r \subset \bg$ such that
$\bg = V_1 \oplus ... \oplus V_r$, with $[V_1,V_i] = V_{i+1}$, $i
= 1,..., r-1$, $[V_1,V_r] = \{0\}$. A natural family of
non-isotropic dilations on $\bg$ associated with this grading is
given by $\Delta_\lambda(\xi) = \lambda \xi_1 + ... + \lambda^r
\xi_r$, if $\xi = \xi_1 + ... + \xi_r \in \bg$. Using the global
diffeomorphism $\exp : \bg \to \bG$, one then lifts these
dilations to the one-parameter family of group automorphisms
\begin{equation}\label{dil}
\delta_\lambda(g)\ =\ \exp \circ \Delta_\lambda \circ
\exp^{-1}(g)\ , \quad\quad\quad \lambda > 0\ . \end{equation}

The homogeneous dimension associated with the dilations
 $\{\delta_\lambda\}_{\lambda>0}$ is given by
 \[
 Q\ =\ \sum_{j = 1}^r j\ dim\ V_j\ .
 \]

Such number often replaces the topological dimension $N =
\sum_{j=1}^r dim\ V_j$ in the analysis of $\bG$, see for instance
Corollary \ref{C:Rellich2} below. Its relevance
 is expressed by the fact that, if $dg$ denotes the bi-invariant
Haar  measure on $\bG$ obtained by pushing forward through the
exponential mapping the Lebesgue measure on $\bg$, then
$d(\delta_\lambda(g)) = \lambda^Q dg$. Here, bi-invariant means with
respect to the operators of left- and right-translation $L_g(g') = g
g'$, $R_g(g') = g' g$, on $\bG$.

The gauge pseudo-distance $\rho(g,g')$ is defined as follows. Let
$|\cdot |$ denote the Euclidean distance to the origin  on
$\mathfrak g$. For $\xi = \xi_1 + \cdots + \xi_r \in \mathfrak g$,
$\xi_i\in V_i$, one lets
\begin{equation}\label{gauge}
|\xi|_{\mathfrak g}\ =\ \left(\sum_{i=1}^r
|\xi_i|^{2r!/i}\right)^{2r!},\quad\quad\quad |g|_{\bG}\ =\
|\exp^{-1}\ g|_{\mathfrak g}, \quad\quad g\in \bG .
\end{equation}
The pseudo-distance on $\bG$ associated to $|\cdot|_{\bG}$ is given
by
\begin{equation}\label{pseudod}
\rho(g, g')\ =\ |g^{-1}\ g'|_{\bG}.
\end{equation}

With a slight abuse of notation when we write $\rho(g)$ we indicate
$\rho(g,e)$, where $e$ is the group identity. Since the function
$g\to \rho(g)$ is homogeneous of degree one with respect to the
non-isotropic dilations \eqref{dil}, we have for the gauge
pseudo-ball,
\begin{equation}\label{balls} B(g,R) = \{g'\in \bG\mid \rho(g',g)<R\} ,
\end{equation}
that $|B(g,R)| = |B(0,1)| R^Q$. Let $m_j = dim\ V_j$, $j=1,...,r$,
and for each $j$ denote by $\{e_{j,s}\}$, $s = 1,...,m_j$, an
orthonormal basis of $V_j$. The vectors $\{e_{j,s}\}$ constitute an
orthonormal basis of $\bg$. Because of the special role played by
the first two layers $V_1, V_2$ in the grading of $\bg$, it is
convenient to have a simpler notation for the elements of their
basis. We thus set henceforth $m = dim\ V_1$, $k = dim\ V_2$, and
will indicate with $\{e_1,...,e_m\}$,
$\{\epsilon_1,...,\epsilon_k\}$ the corresponding basis of $V_1,
V_2$. Denoting with $(L_g)_*$ the differential of left-translations,
we define a family of left-invariant vector fields on $\bG$ by
letting
\begin{equation}\label{vf}
\begin{cases}
X_i(g)\ =\ (L_g)_*(e_i)\ , \ i=1,...,m\ ,\ \ T_s(g)\ =\ (L_g)_*(\epsilon_s)\ , \ s=1,...,k\ , \\
X_{j,s}(g)\ =\ (L_g)_*(e_{j,s})\ ,\ j = 3,..., r\ , \ s = 1,...,m_j\
.
\end{cases}
\end{equation}

Hereafter, we assume that $\bG$ is endowed with a left-invariant
Riemannian metric $<\cdot,\cdot>$ with respect to which the vector
fields $X_1,...,X_m, T_1,...,T_k,X_{j,s}$, $j=3,...r$, $s = 1,...,
m_j$, are orthonormal. In view of the grading assumption on $\bg$,
it is clear that the vector fields $X_1,...,X_m$ generate the Lie
algebra of all left-invariant vector fields on $\bG$. They generate
a sub-bundle $H\bG$ of the tangent bundle $T\bG$ which is usually
called the horizontal bundle.

We next recall some basic concepts from the sub-Riemannian geometry
of an hypersurface in a Carnot group $\bG$. For a detailed account
we refer the reader to \cite{DGN2}. We consider the Riemannian
manifold $M = \bG$ with the left-invariant metric tensor with
respect to which $X_1,...,X_m, ... , X_{r,m_r}$ is an orthonormal
basis, the corresponding Levi-Civita connection $\nabla$ on $\bG$,
and the horizontal Levi-Civita connection $\nabla^H$. Let $\Om
\subset \bG$ be a bounded $C^k$ domain, with $k\geq 2$. We denote by
$\n$ the Riemannian outer normal to $\p \Om$, and define the
so-called \emph{angle function} on $\p \Om$ as follows
\begin{equation}\label{ps}
W\ =\ |\bN^H|\ =\ \sqrt{\sum_{j=1}^m <\n, X_j>^2}\ .
\end{equation}

The \emph{characteristic set} of $\Om$, hereafter denoted by
$\Sigma$, is the compact subset of $\p \Om$ where the continuous
function $W$ vanishes
\begin{equation}\label{csaf}
\Sigma\ =\ \{g\in \p \Om\mid W(g) = 0\}\ .
\end{equation}

The next definition plays a basic role in sub-Riemannian geometry.

\begin{dfn}\label{D:HGauss}
We define the outer \emph{horizontal normal} on $\p \Om$ as follows
\begin{equation}\label{up1}
\bN^H\ =\ \sum_{j=1}^m <\n,X_j> X_j\  ,
\end{equation}
so that $W = |\bN^H|$. The \emph{horizontal Gauss map} $\nuX$ on $\p
\Om$ is defined by
\begin{equation}\label{up2}
\nuX\ =\ \frac{\bN^H}{|\bN^H|}\ ,\quad\quad\quad\quad \text{on}\quad
\p \Om \setminus \Sigma\ .
\end{equation}

\end{dfn}

We note that $\Up$ is the projection of the Riemannian Gauss map on
$\p \Om$ onto the horizontal subbundle $H\bG \subset T \bG$. Such
projection vanishes only at characteristic points, and this is why
the horizontal Gauss map is not defined on $\Sigma$. The following
definition is taken from \cite{DGN2}, but the reader should also see
\cite{HP} for a related notion in the more general setting of
vertically rigid spaces.

\begin{dfn}\label{D:HMC}
The \emph{horizontal} or $H$-\emph{mean curvature} of $\p \Om$ at a
point $g_0\in \p \Om \setminus \Sigma$ is defined as
\[
\mathcal H\ =\ \sum_{i=1}^{m-1} <\nabla^H_{\boldsymbol e_i}
\boldsymbol e_i,\nuX>\ ,
\]
where $\{\boldsymbol e_1,...,\boldsymbol e_{m-1}\}$ denotes an
orthonormal basis of the horizontal tangent bundle $T_H \p \Om
\overset{def}{=} T \p \Om \cap H\bG$ on $\p \Om$. If instead $g_0\in
\Sigma$, then we define $\mathcal H(g_0) = \underset{g\to g_0}{\lim}
\mathcal H(g)$, provided that the limit exists and is finite.
\end{dfn}

We next consider the following nonlinear operator
\begin{equation}\label{infty}
\Delta_{H,\infty} u\ \overset{def}{=}\ \sum_{i,j=1}^m u_{,ij}\ X_i
u\ X_j u\ ,
\end{equation}
which by analogy with its by now classical Euclidean ancestor we
call the \emph{horizontal $\infty$-Laplacian}. This operator has
been recently studied by various people, see e.g. \cite{Bi},
\cite{BiC}, \cite{Wa3}, \cite{GT}. The reason for introducing the
operator $\Delta_{H,\infty}$ is in the following proposition which
is often useful in computing the $H$-mean curvature. To state it we
recall the notion of a defining function for $\Om$. We
 consider a $C^2$ bounded open set $\Om \subset \bG$ and we
 assume for convenience that there exists a globally defined $\phi \in C^2(\bG)$ (a defining function) such that
\begin{equation}\label{defining}
\Om\ =\ \{g\in \bG\ \mid\ \phi(g)\ <\ 0\}\ ,
\end{equation}
and for which $|\nabla \phi| \geq \alpha > 0$ in an open
neighborhood $\mathcal O$ of $\p \Om$, where $\nabla \phi$ denotes
the Riemannian gradient of $\phi$. The Riemannian outer unit normal
to $\p \Om$ is presently given by $\n = \nabla \phi/|\nabla \phi|$.
We observe that
\begin{equation}\label{up3}
|\Up|\ =\ \frac{|\nh\phi|}{|\nabla \phi|}\ ,
\end{equation}
and that on $\p \Om \setminus \Sigma$ one has
\begin{equation}\label{nuX}
\nuX\ =\ \frac{\nh\phi}{|\nh\phi|}\ .
\end{equation}

The next result is Proposition 9.12 in \cite{DGN2}.

\begin{prop}\label{P:XMCG}
At every point of $\p \Om\setminus \Sigma$ one has in terms of a
local defining function $\phi$ of $S$
\begin{align*}
|\nh \phi|^3\ \mathcal H\ =\  |\nh\phi|^2\ \sul \phi\ -\
\Delta_{H,\infty} \phi\ .
\end{align*}
\end{prop}

We will need the following lemma.

\begin{lemma}\label{L:convlevelset}
Let $u\in C^2(\bG)$ be $H_2$-convex, then for every $s\in \R$ such
that the level set
\[
\Om_s\ =\ \{g\in \bG \mid u(g) < s\}
\] is a $C^2$ domain, the $H$-mean curvature of
$\mathcal E_s = \p \Om_s$ (wherever it is defined) is nonnegative.
\end{lemma}

\begin{proof}[\textbf{Proof}]
Recall that the hypothesis that $u$ be $H_2$-convex means that $\sul
u \geq 0$, and that moreover \begin{equation}\label{h2}
 (\sul u)^2\
-\ ||\nabla^2_H u||^2\ \geq\ 0\ .
\end{equation}

According to Proposition \ref{P:XMCG} it suffices to show that on
$\mathcal E_s$ one has $|\nh u|^2 \sul u - \Delta_\infty u\ \geq\
0$. On the other hand, Schwarz inequality gives
\[
\Delta_{H,\infty} u\ =\ \sij u_{,ij} X_i u X_j u\ \leq\ ||\nabla^2_H
u||\ |\nh u|^2\ \leq\ \sul u\ |\nh u|^2\ ,
\]
where in the last inequality we have used \eqref{h2}.

\end{proof}

Given an open set $\Om \subset \bG$ denote by
\[
\mathcal F(\Om)\ =\ \{\phi = \sum_{j=1}^m \phi_j X_j \in
C^1_0(\Om,H\bG) \mid ||\phi||_\infty = \underset{g\in \Om}{\sup}
(\sum_{j=1}^m \phi_j^2)^{1/2} \leq 1\}\ .
\]

The $H$-perimeter of a measurable set $E\subset \bG$ with respect to
$\Om$ was defined in \cite{CDG} as
\[
P_H(E;\Om) \ =\ \underset{\phi\in \mathcal F(\Om)}{\sup} \int_{E\cap
\Om} div_H \phi\ dg\ .
\]

If $E$ is a bounded open set of class $C^1$, then the divergence
theorem gives
\[
P_H(E;\Om) = \underset{\phi\in \mathcal F(\Om)}{\sup} \int_{\p E\cap
\Om} <\n,X_j> \phi_j\ d\sigma = \underset{\phi\in \mathcal
F(\Om)}{\sup} \int_{\p E\cap \Om} <\bN^H,\phi>  d\sigma = \int_{\p
E\cap \Om} |\bN^H| d\sigma\ ,
\]
where $d\sigma$ is the Riemannian surface measure on $\p E$. It is
clear form this formula that the measure on $\p E$, defined by
\[ \sigma_H(\p E \cap \Om)\ \overset{def}{=}\ P_H(E;\Om)\  \]
on the open sets of $\p E$, is absolutely continuous with respect to
$\sigma$, and its density is represented by the angle function $W$
of $\p E$. We formalize this observation in the following
definition.

\begin{dfn}\label{D:permeas}
Given a bounded domain $E \subset \bG$ of class $C^1$, with angle
function $W$ as in \eqref{csaf}, we will denote by
\begin{equation}\label{permeasure}
d \sigma_H\ =\ |\Up|\ d\sigma\ =\ W\ d\sigma\ ,
\end{equation}
the $H$-perimeter measure supported on $\p E$.
\end{dfn}

\section{\textbf{A Bochner type identity}}\label{S:bochner}

In this section we establish a sub-Riemannian version of the
classical Bochner identity \eqref{bochner} for the sub-Laplacian of
the square of the length of the horizontal gradient of a function on
a Carnot group $\bG$, see Proposition \ref{P:B}. A deeper CR version
of such formula first appeared in the beautiful paper by A.
Greenleaf \cite{Gre}. We begin by finding the formula which
expresses the connection between the Hilbert-Schmidt norm of the
symmetrized, and that of the un-symmetrized horizontal Hessian.

\begin{lemma}\label{L:symunsym}
Let $u\in \Gamma^2(\bG)$, then one has
\[
\sum_{i,j=1}^m (X_i X_j u)^2\ =\ ||\nabla_H^2 u||^2\ +\ \frac{1}{4}\
\sum_{i,j=1}^m ([X_i,X_j]u)^2 \ .
\]
\end{lemma}

\begin{proof}[\textbf{Proof}]
Notice that the un-symmetrized and the symmetrized second
derivatives are connected by the formula
\begin{equation}\label{symunsym} X_iX_j u\ =\ u_{,ij}\ +\
\frac{1}{2}\ [X_i,X_j] u\ .
\end{equation}

We obtain from \eqref{symunsym}
\[
\sum_{i,j=1}^m (X_i X_j u)^2\ =\ ||\nabla_H^2 u||^2\ +\ \frac{1}{4}\
\sum_{i,j=1}^m ([X_i,X_j]u)^2 \ +\ \sij u_{,ij}\ [X_i , X_j]u\ .
\]

To reach the conclusion, it is now enough to observe that, thanks to
the skew-symmetry of the matrix $\{[X_i,X_j]\}$, we have
\[ \sij u_{,ij}\ [X_i , X_j]u\ =\ \sum_{i<j} u_{,ij}\ [X_i , X_j]u\
+\ \sum_{i>j} u_{,ij}\ [X_i , X_j]u\ =\ 0\ .
\]

\end{proof}

\begin{lemma}\label{L:comm}
For $u\in \Gamma^2(\bG)$ one has
\[
 \sij\  X_i X_j u\ [X_i,X_j] u\ =\ \frac{1}{2}\ \sij ([X_i,X_j]
u)^2\ .
\]
\end{lemma}

\begin{proof}[\textbf{Proof}]
To check this formula we proceed as follows
\begin{align}\label{g4} & \sij\  X_i X_j u\ [X_i,X_j] u\ =\
\sum_{i<j} X_i X_j u\ [X_i,X_j] u\ +\ \sum_{i>j} X_i X_j u\
[X_i,X_j] u
\\
& =\ \sum_{i<j} X_i X_j u\ [X_i,X_j] u\ -\ \sum_{i<j} X_j X_i u\
[X_i,X_j] u \notag\\
& =\ \sum_{i<j} ([X_i,X_j] u)^2\ =\ \frac{1}{2}\ \sij ([X_i,X_j]
u)^2\ . \notag
\end{align}

\end{proof}

\begin{prop}\label{P:B}
Let $\bG$ be a Carnot group, $u\in \Gamma^3(\bG)$, then the
following sub-Riemannian Bochner formula holds
\begin{align*}
 \frac{1}{2}\ \sul (|\nh u|^2)\ & =\ <\nh u , \nh(\sul)>\ +\
||\nabla^2_H u||^2\ +\ \frac{1}{4}\ \sum_{i,j=1}^m ([X_i,X_j]u)^2
\\
&  +\ 2\ \sij X_ju\ [X_i,X_j] X_iu \ +\ \sij X_j u\ [X_i,[X_i,X_j]]
u\ .\notag
\end{align*}
When $\bG$ is of step $2$, then for every $i,j =1,...,m,$ one has
$[X_i,[X_i,X_j]] = 0$ in the last term in the right-hand side of the
above identity. In particular, when $\bG = \Hn$, then we have
\begin{align}\label{B2}
 \frac{1}{2}\ \sul (|\nh u|^2)\ & =\ <\nh u , \nh(\sul u)>\ +\
||\nabla_H^2 u||^2\ +\ \frac{3}{2}\ n\ (T u)^2
\\
&  +\ 2\ \sij X_j u [X_i,X_j] X_iu\ .\notag
\end{align}
\end{prop}

\begin{proof}[\textbf{Proof}]
We observe that for any function $F$ we have
\begin{equation}\label{square}
\sul (F^2)\ =\ 2\ F\ \sul F\ +\ 2\ |\nh F|^2\ .
\end{equation}

Applying \eqref{square} to $F = X_ju$ we obtain
\begin{equation}\label{B2}
\frac{1}{2}\ \sul (|\nh u|^2)\ =\ \frac{1}{2}\ \sum_{j=1}^m \sul
((X_j u)^2)\ =\ \sum_{j=1}^m X_j u\ \sul(X_ju)\ +\ \sum_{i,j=1}^m
(X_iX_ju)^2\ .
\end{equation}

We next compute $\sul(X_ju)$. One has
\begin{align}\label{B3}
\sul(X_ju)\ & =\ \sum_{i=1}^m X_i X_i X_ju\ =\ \sum_{i=1}^m X_i
(X_j X_iu\ +\ [X_i,X_j]u)
\\
& =\ \sum_{i=1}^m (X_j X_i \ +\ [X_i,X_j])X_iu\ +\ \sum_{i=1}^m
X_i [X_i,X_j]u
\notag\\
& =\ X_j(\sul u)\ +\ \sum_{i=1}^m [[X_i,X_j], X_i]u\ +\ 2\
\sum_{i=1}^m [X_i,X_j]X_i u\ . \notag
\end{align}

On the other hand, Lemma \ref{L:symunsym} gives
\begin{equation}\label{B4bis}
\sum_{i,j=1}^m (X_i X_j u)^2\ =\ ||\nabla_H^2 u||^2\ +\ \frac{1}{4}\
\sum_{i,j=1}^m ([X_i,X_j]u)^2 \ .
\end{equation}

Substituting \eqref{B3}, \eqref{B4bis} in \eqref{B2} we find
\begin{align*}
 \frac{1}{2}\ \sul (|\nh u|^2)\ & =\ <\nh u , \nh(\sul)>\ +\
||\nabla^2_H u||^2\ +\ \frac{1}{4}\ \sum_{i,j=1}^m ([X_i,X_j]u)^2
\\
&  +\ 2\ \sij X_ju\ [X_i,X_j] X_iu \ +\ \sij X_j u\ [X_i,[X_i,X_j]]
u\ ,\notag
\end{align*}
which gives the desired conclusion.

\end{proof}

%>>>>>>>>>>>>>>>>>>>>>>>>>>>>>>>>>>>>>>>>>>>>>

\section{\textbf{Geometric second derivative estimates}}\label{S:2nd}

In this section using the horizontal Bochner identity in Proposition
\ref{P:B} we prove the various results stated in the introduction.

\begin{proof}[\textbf{Proof of Theorem \ref{T:grisvard}}]
We begin by observing that, if we denote by $\n$ the outer unit
Riemannian normal on $\p \Om$, then the assumptions $u \leq 0$ in
$\Om$ and $u=0$ on $\p \Om$ imply
\begin{equation}\label{nu}
\nabla u\ =\ |\nabla u|\ \n\ ,\quad\quad\quad \text{on}\quad \p \Om\
.
\end{equation}

Next, we rewrite the identity in Proposition \ref{P:B} as follows
\begin{align*}
 \frac{1}{2}\ \sul (|\nh u|^2)\ & =\ <\nh u , \nh(\sul)>\ +\
||\nabla^2_H u||^2\ +\ \frac{1}{4}\ \sum_{i,j=1}^m ([X_i,X_j]u)^2
\\
&  +\ 2\ \sij X_ju\ X_i [X_i,X_j] u \ +\ 2\ \sij X_ju\
[[X_i,X_j],X_i]u \notag\\
& +\ \sij X_j u\ [X_i,[X_i,X_j]] u\ . \notag
\end{align*}

This gives
\begin{align}\label{Bmod}
\frac{1}{2}\ \sul (|\nh u|^2)\ &  =\ <\nh u , \nh(\sul)>\ +\
||\nabla^2_H u||^2\ +\ \frac{1}{4}\ \sum_{i,j=1}^m ([X_i,X_j]u)^2
\\
&  +\ 2\ \sij X_ju\ X_i [X_i,X_j] u \ +\  \sij X_ju\
[[X_i,X_j],X_i]u\ . \notag
\end{align}

We now integrate the identity \eqref{Bmod} on $\Om$
\begin{align}\label{g1}
& \frac{1}{2}\ \int_\Om \Delta_H(|\nh u|^2)\ dg\ =\ \int_\Om <\nh
u , \nh (\Delta_H u)>  dg
\\
 & +\ \int_\Om
||\nabla^2_H u||^2\ dg\ +\ \frac{1}{4}\ \sum_{i,j=1}^m \int_\Om
([X_i , X_j] u)^2\  dg \notag\\
&  +\ 2\ \sij \int_\Om X_ju\ X_i [X_i,X_j] u \ dg\ +\  \sij \int_\Om
X_ju\ [[X_i,X_j],X_i]u\ dg\ . \notag
\end{align}

Using \eqref{nu} we have from the divergence theorem
\begin{equation}\label{g2}
\frac{1}{2}\ \int_\Om \Delta_H(|\nh u|^2)\ dg\ =\ \frac{1}{2}\
\int_{\p \Om} \frac{<\nh (|\nh u|^2), \nh u>}{|\nabla u|}\ d\sigma\
.
\end{equation}

Also, again from \eqref{nu}, we find
\begin{equation}\label{g3}
\int_\Om <\nh u , \nh (\Delta_H u)>  dg\ =\ \int_{\p \Om} \frac{|\nh
u|^2 \Delta_H u}{|\nabla u|}\ d\sigma \ -\ \int_\Om (\Delta_H u)^2\
dg\ .
\end{equation}

Finally, we have
\begin{align}\label{g4}
& 2\ \sij \int_\Om X_j u \ X_i[X_i,X_j] u\ dg\  =\ 2\ \sij \int_{\p
\Om} \frac{[X_i,X_j]u\ X_i u X_j u}{|\nabla u|}\ d\sigma
\\
& - \ 2\ \sij \int_\Om [X_i,X_j] u\ X_i X_ju\ dg\ =\ -\ 2\ \sij
\int_\Om [X_i,X_j] u\ X_i X_ju\ dg\ \notag\ ,
\end{align}
where to eliminate the boundary integral we have used the
skew-symmetry of the matrix $\{[X_i,X_j]u\}_{i,j=1,...,m}$.
Substituting \eqref{g2}-\eqref{g4} into \eqref{g1} we conclude
\begin{align*}
& \int_\Om \bigg\{(\Delta_H u)^2 -  ||\nabla^2_H u||^2\bigg\} dg\
 +\ \frac{3}{4}\ \sum_{i,j=1}^m \int_\Om ([X_i , X_j] u)^2\ dg\ +\ \sij \int_\Om
X_iu\ [[X_i,X_j],X_j]u\ dg\
\\
& =\ \int_{\p \Om} \frac{\big\{|\nh u|^2 \sul u - \Delta_{H,\infty}
u\big\}}{|\nabla u|} d\sigma\ =\ \int_{\p \Om} \mathcal H\ |\nh u|^2
d\sigma_H\ ,
\end{align*}
where in the last equality we have used Proposition \ref{P:XMCG} and
Definition \ref{D:permeas}. This gives the desired conclusion.

\end{proof}

\begin{cor}\label{C:grisvard}
Let $\bG$ be a Carnot group of step $r=2$, and consider a $C^2$
bounded open set $\Om\subset \bG$. Let $u\in \Gamma^2(\overline
\Om)$ with $u\leq 0$ in $\Om$ and $u = 0$ on $\p \Om$. One has
\begin{align}\label{cor1}
& \int_\Om ||\nabla^2_H u||^2\ dg\ =\ \int_\Om (\Delta_H u)^2\ dg\
+\ \frac{3}{4} \sij \int_\Om ([X_i,X_j]u)^2\ dg
\\
&   -\ \int_{\p \Om} \mathcal H\ |\nabla_H u|^2\ d\sigma_H\
.\notag
\end{align}
If instead $u\in \Gamma^2_0(\Om)$, then one obtains
\begin{align}\label{cor2}
\int_\Om ||\nabla^2_H u||^2\ dg\ & =\ \int_\Om (\Delta_H u)^2\ dg\
+\ \frac{3}{4} \sij \int_\Om ([X_i,X_j]u)^2\ dg\ .
\end{align}
\end{cor}

We can now present the

\begin{proof}[\textbf{Proof of Theorem \ref{T:commH}}]
We notice that, since $u$ is $H_2$-convex, then $\Delta_H u\geq 0$,
and $(\Delta_H u)^2 - ||\nabla^2_H u||^2 \geq 0$ in $\Om$. In
particular, since by assumption $u= 0$ on $\p \Om$, from Bony's weak
maximum principle \cite{Bo} we infer that $u \leq 0$ in $\Om$, and
therefore we can apply Theorem \ref{T:grisvard}. The desired
conclusion now follows from \eqref{1} in Theorem \ref{T:grisvard}.

\end{proof}

We next want to control the commutator term in the right-hand side
of \eqref{cor1} in Corollary \ref{C:grisvard}. To reach this goal we
will make use of a sub-Riemannian Rellich identity discovered in
\cite{GV}. In the following results, $\Om$ will indicate a piecewise
$C^1$ bounded open subset of a Carnot group $\bG$ with outer unit
normal $\n$ and surface measure $\sigma$.

\begin{thrm}\label{T:Rellich}
For $u\in \Gamma^{2}(\overline{\Om})$ one has
\begin{align*}
&  2\ \int_{\partial{\Om}}\ \zeta u\ <\nh u, \bN^H> d\sigma
 + \int_\Om \ div_{\bG} \zeta\ |\nh u|^2\ dg
\\
&  -\ 2\ \sum_{i=1}^m\ \int_\Om\ X_iu\ [X_i,\zeta]u\ dg \
 -\ 2\ \int_\Om\ \zeta u\ \Delta_H u\ dg
\notag\\
& =\ \int_{\partial{\Om}}\ |\nh u|^2\ <\zeta,\boldsymbol \nu>\
d\sigma , \notag
\end{align*}
where $\zeta$ is a $C^1$ vector field on $\bG$.
\end{thrm}

\begin{cor}\label{C:Rellich}
Let $u\in \Gamma^{2}(\overline{\Om})$ and assume, in addition, that
$u = 0$ on $\p \Om$. One has
\begin{align*}
&   \int_{\partial{\Om}}  |\nh u|^2 <\zeta , \boldsymbol \nu>
d\sigma
 + \int_\Om \ div_{\bG} \zeta\ |\nh u|^2\ dg
 \\
&  -\ 2\ \sum_{i=1}^m\ \int_\Om\ X_iu\ [X_i,\zeta]u\ dg \
 -\ 2\ \int_\Om\ \zeta u\ \Delta_H u\ dg\ =\ 0\ .
\notag
\end{align*}
\end{cor}

In what follows we indicate with $\mathcal Z$ the infinitesimal
generator of the non-isotropic dilations \eqref{dil}. We note that
in the exponential coordinates it is given by
\begin{equation}\label{Z}
\mathcal Z\ =\ \sum_{i=1}^m\ x_i(g)\ X_i\ +\ 2\ \sum_{s=1}^k\
t_s(g) \ T_s\ +\  \sum_{j=3}^r j\ \sum_{s=1}^{m_j} x_{j,s}(g)\
X_{j,s}\ .
\end{equation}

When the step of the group is $r=2$ the third sum in the right-hand
side of \eqref{Z} does not appear.
\begin{equation}\label{divZ}
[X_i,\mathcal Z]\ =\ X_i\ ,\ i=1,...,m\ ,\ \ \ div_\bG\ \mathcal Z\
=\ Q\ .
\end{equation}

For a proof of the first identity in \eqref{divZ} see Lemma 2.1 in
\cite{DG}. The second identity follows by using the expression
\eqref{Z} of $\mathcal Z$ in the exponential coordinates. Choosing
$\eta = \mathcal Z$ in Corollary \ref{C:Rellich}, and using
\eqref{divZ} we easily obtain.

\begin{cor}\label{C:Rellich2}
Let $u\in \Gamma^{2}(\overline{\Om})$ and assume, in addition, that
$u = 0$ on $\p \Om$. One has
\[
 \int_{\partial{\Om}}  |\nh u|^2 <\mathcal Z , \boldsymbol \nu> d\sigma\
 +\ (Q - 2) \int_\Om  |\nh u|^2\ dg\ =\ 2 \int_\Om\ \mathcal Z u\ \Delta_H u\ dg\  .
\]
\end{cor}

Using Corollary \ref{C:Rellich2} we can now prove the following
useful estimate.

\begin{lemma}\label{L:rellich3}
Suppose that $\bG$ be a Carnot group of step $r=2$. Under the
hypothesis of Corollary \ref{C:Rellich2} on the function $u$ there
exists a constant $C = C(\bG,\Om)>0$ such that for any $\epsilon >0$
one has
\begin{align*}
& \int_{\partial{\Om}}  |\nh u|^2 <\mathcal Z , \boldsymbol \nu>
d\sigma\
 +\ \frac{(Q - 2)}{2} \int_\Om  |\nh u|^2\ dg
 \\
 & \leq\ C \left\{\left(1 +
 \frac{1}{\epsilon}\right) \int_\Om (\Delta_Hu)^2\ dg\ +\ \epsilon
 \sum_{i,j=1}^m \int_\Om ([X_i,X_j]u)^2\ dg\    \right\}\ .
\end{align*}
\end{lemma}

\begin{proof}[\textbf{Proof}]
Since $\bG$ has step $r=2$, from the bracket generating assumption
for every $s=1,...,k$ there exist $\alpha_{i,j}^s\in \R$,
$i,j=1,...,m$, such that
\[
T_s\ =\ \sum_{i,j = 1}^m \alpha_{i,j}^s [X_i,X_j]\ .
\]

Therefore, it is possible to find $\beta_s>0$ such that
\[
|T_su|\ \leq\ \beta_s \left(\sum_{i,j=1}^m
([X_i,X_j]u)^2\right)^{1/2}\ .
\]

From this estimate, from \eqref{Z} and from the boundedness of $\Om$
we conclude that there exists $C = C(\bG,\Om)>0$ such that one has
\[
|\mathcal Zu\ \Delta_H u|\ \leq\ C\ \left\{|\nh u|\ +\
\left(\sum_{i,j=1}^m ([X_i,X_j]u)^2\right)^{1/2}\right\}\ |\Delta_H
u|\ ,\quad\quad\quad \text{in}\ \overline \Om\ .
\]

Inserting this estimate in the identity of Corollary
\ref{C:Rellich2}, for every $\delta, \epsilon >0$ we find
\begin{align*}
& \int_{\partial{\Om}}  |\nh u|^2 <\mathcal Z , \boldsymbol \nu>
d\sigma\
 +\ (Q - 2) \int_\Om  |\nh u|^2\ dg\ \leq\ 2 \int_\Om\ |\mathcal Z u\ \Delta_H u|\ dg
 \\
 & \leq\ C \delta \int_\Om |\nh u|^2\ dg\ +\ \frac{C}{\delta}\ \int_\Om (\Delta_Hu)^2\ dg
 \\
 & +\  C \epsilon  \sum_{i,j=1}^m \int_\Om
([X_i,X_j]u)^2 \ dg\ +\ \frac{C}{\epsilon}\ \int_\Om (\Delta_Hu)^2\
dg\ .
\end{align*}

Choosing now $\delta>0$ such that $C \delta = \frac{Q-2}{2}$ we
obtain the desired conclusion (with a possibly different constant $C
= C(\bG,\Om)>0$).

\end{proof}

We can now provide the

\begin{proof}[\textbf{Proof of Theorem \ref{T:est}}]
We start with the inequality \eqref{commH1} in Theorem
\ref{T:commH}, which gives
\begin{align}\label{finalrush}
\frac{3}{4} \sij \int_\Om ([X_i,X_j]u)^2\ dg\ & \leq\ \int_{\p \Om}
\mathcal H\ |\nabla_H u|^2\ d\sigma_H\ \leq\ M\ \int_{\p \Om} |\nh
u|^2 d\sigma_H\ ,
\end{align}
where we have used the hypothesis \eqref{boundedmc}. According to
Lemma \ref{L:rellich3} if $\Om$ satisfies the hypothesis
\eqref{starlike} we obtain for any $\epsilon >0$
\begin{align*}
& \alpha\ \int_{\partial{\Om}}  |\nh u|^2\ W\ d\sigma\
 +\ \frac{(Q - 2)}{2} \int_\Om  |\nh u|^2\ dg
\\
& \leq\ \int_{\partial{\Om}}  |\nh u|^2 <\mathcal Z , \boldsymbol
\nu> d\sigma\
 +\ \frac{(Q - 2)}{2} \int_\Om  |\nh u|^2\ dg
 \\
 & \leq\ C \left\{\left(1 +
 \frac{1}{\epsilon}\right) \int_\Om (\Delta_Hu)^2\ dg\ +\ \epsilon
 \sum_{i,j=1}^m \int_\Om ([X_i,X_j]u)^2\ dg\    \right\}\ .
\end{align*}

Keeping \eqref{permeasure} in mind, we have proved that for any
$\epsilon>0$
\[
\int_{\partial{\Om}}  |\nh u|^2\ d\sigma_H\ \leq\ C(\bG,\alpha)\
\left\{\left(1 +
 \frac{1}{\epsilon}\right) \int_\Om (\Delta_Hu)^2\ dg\ +\ \epsilon
 \sum_{i,j=1}^m \int_\Om ([X_i,X_j]u)^2\ dg\    \right\}\ .
\]

Combining this estimate with \eqref{finalrush} we obtain
\begin{align}\label{finalrush2}
\sij \int_\Om ([X_i,X_j]u)^2\ dg\ \leq\ C(\bG,M,\alpha)\
\left\{\left(1 +
 \frac{1}{\epsilon}\right) \int_\Om (\Delta_Hu)^2\ dg\ +\ \epsilon
 \sum_{i,j=1}^m \int_\Om ([X_i,X_j]u)^2\ dg\    \right\}\ .
\end{align}

Choosing $\epsilon>0$ in \eqref{finalrush2} such that $\epsilon
C(\bG,M,\alpha)<1$ we finally reach the desired conclusion.

\end{proof}

We close with a proposition which provides a significant class of
domains which satisfy the two geometric hypothesis in Theorem
\ref{T:est}.

\begin{prop}\label{P:gaugeball}
In a Carnot group $\bG$ of arbitrary step consider a gauge
pseudo-ball $B_R = \{g\in \bG\mid \rho(g)<R\}$, where $ \rho$ is the
Folland-Stein gauge \eqref{gauge}, \eqref{pseudod}. There exists $C
= C(\bG)>0$, $\alpha = \alpha(\bG)>0$ such that
\begin{equation}\label{boundedmcg}
\underset{\p B_R}{\sup}\ |\mathcal H|\ \leq\ \frac{C}{R}\ ,
\end{equation}
and
\begin{equation}\label{starlikeg}
\underset{\p B_R}{\inf}\ <\mathcal Z,\n>\ \geq\ \alpha\ R\ W\ .
\end{equation}
\end{prop}

\begin{proof}[\textbf{Proof}]
The outer unit normal to $B_R$ at a point of its boundary is given
by $\n = \frac{\nabla \rho}{|\nabla \rho|}$. Since the function
$\rho$ is homogeneous of degree one with respect to the
non-isotropic group dilations, from the Euler type formula for
Carnot groups we obtain on $\p B_R$
\[
<\mathcal Z,\n>\ =\ <\mathcal Z,\frac{\nabla \rho}{|\nabla \rho|}>\
=\ \frac{\mathcal Z\rho}{|\nabla \rho|}\ =\ \frac{\rho}{|\nabla
\rho|}\ =\ \frac{R}{|\nabla \rho|}\ .
\]

On the other hand, since $\rho\in C^\infty(\bG\setminus\{e\})$, and
since $|\nabla^H \rho|$ is homogeneous of degree zero, we have for
every $g\not= e$
\[
W(g)\ =\ \frac{|\nabla^H \rho(g)|}{|\nabla \rho(g)|}\ \leq\
\frac{\underset{\rho(g')=1}{sup} |\nabla^H \rho(g')|}{|\nabla
\rho(g)|}\ =\ \frac{C(\bG)}{|\nabla \rho(g)|}\ .
\]

We thus obtain on $\p B_R$
\[ <\mathcal Z,\n>\ \geq\ C(\bG)^{-1} R\
W\ =\ \alpha\ R\ W\ .
\]

This proves \eqref{starlikeg}. To prove the qualitative estimate
\eqref{boundedmcg} we again employ homogeneity considerations.
According to Proposition \ref{P:XMCG} we have on $\p B_R$
\begin{align*}
|\nh \rho|^3\ \mathcal H\ =\  |\nh\rho|^2\ \sul \rho\ -\
\Delta_{H,\infty} \rho\ .
\end{align*}

Now, $\sul \rho$ and $\Delta_{H,\infty}\rho$ both have homogeneity
$-1$, and hence so does $\mathcal H$. We thus find on $\p B_R$
\[
|\mathcal H(g)|\ \leq\ \frac{1}{\rho(g)}\ \underset{\rho(g')=1}{sup}
|\mathcal H(g')|\ =\ \frac{C(\bG)}{R}\ ,
\]
which establishes \eqref{boundedmcg}.

\end{proof}

\end{document}